\documentclass[11pt,notitlepage,twoside,a4paper]{amsart}
 \usepackage{amsfonts}

\usepackage{amsmath,amssymb,enumerate}

\usepackage{epsfig,fancyhdr,color}

\usepackage{amssymb}
\usepackage{amsmath,amsthm}
\usepackage{latexsym}
\usepackage{amscd}
\usepackage{psfrag}
\usepackage{graphicx}
\usepackage[latin1]{inputenc}
\usepackage[all]{xy}
\usepackage[mathcal]{eucal}

\definecolor{NoteColor}{rgb}{1,0,0}

\renewcommand{\textsc}{\textcolor{red}}

\newtheorem*{theorem 1}{\rm\bf Proposition 1}
\newtheorem*{theorem 2}{\rm\bf Proposition 2}

\theoremstyle{definition}

\theoremstyle{remark}

\def\interieur#1{\mathord{\mathop{\kern 0pt #1}\limits^\circ}}

\title[Commentary]
{A commentary on Teichm\"uller's paper \\Ver\"anderliche Riemannsche Fl\"achen\\ (Variable Riemann Surfaces)\\
Deutsche Math. 7, 344-359 (1944)} 

\author{A. A'Campo-Neuen}
\address{Annette A'Campo-Neuen, Mathematisches Institut,
Universit\"at Basel,
Rheinsprung 2,
CH-4051 Basel,
Switzerland} 
\email{Annette.Acampo@unibas.ch}

\author{N. A'Campo}
\address{Norbert A'Campo, Mathematisches Institut,
Universit\"at Basel,
Rheinsprung 2,
CH-4051 Basel,
Switzerland} 
\email{Norbert.ACampo@unibas.ch}

\author{L. Ji}
\address{Lizhen Ji, University of Michigan, Department of Mathematics,  2074 East Hall, 530 Church Street, 
Ann Arbor, MI 48109-1043, USA} 
\email{lji@umich.edu}

\author{A. Papadopoulos}
\address{Athanase Papadopoulos,  Universit{\'e} de Strasbourg and CNRS, Institut de Recherche Math\'ematique Avanc\'ee,
7 rue Ren\'e Descartes,
 67084 Strasbourg Cedex, France} 
\email{athanase.papadopoulos@math.unistra.fr}

\date{\today}

\begin{document}

 \maketitle



This paper is the last one that Teichm\"uller wrote on the problem of moduli. At most places the paper contains ideas and no technical details. The author presents a completely new approach to Teichm\"uller space, compared to the approach he took in his first seminal paper ``Extremale quasikonforme Abbildungen und quadratische Differentiale" \cite{T1939} (and its sequel \cite{T1943} in which he completed some of the the results stated in \cite{T1939}).

 In the paper \cite{T1939}, Teichm\"uller led the foundations of what we call today Teichm\"uller theory (but without the complex structure), defining its metric and introducing in that theory the techniques of quasiconformal mappings and of quadratic differentials as essential tools. In the present paper, the approach is more abstract, through complex analytic geometry. Teichm\"uller space, equipped with its complex-analytic structure, is characterized here by a certain universal property. A fibre bundle is constructed, an object which today bears the name ``Teichm\"uller universal curve". 

 It seems that Teichm\"uller considered that the methods of his 1939 paper \cite{T1939} could not lead to the definition of a complex structure for Teichm\"uller space.\footnote{One can say \emph{a posteriori} that Teichm\"uller was wrong in that, since we know now that the complex structure of Teichm\"uller space can be defined using quasiconformal mappings. Let us quote in this respect Ahlfors, who was the first to derive the complex structure of Teichm\"uller space from the quasiconformal theory. Ahlfors writes in \cite{Ahlfors1960}: ``Teichm\"uller states explicitly, in 1944, that his metrization is of no use for the construction of an analytic structure. The present author disagrees and will show that the metrization and the corresponding parametrization are at least very convenient tools for setting up the desired structure." We also add by the way that we know today that conversely, the Teichm\"uller metric can be recovered from the complex structure, since Royden showed that the Teichm\"uller metric coincides with the Kobayashi metric. Today, the complex structure of Teichm\"uller space is usually presented using the variational theory of the Beltrami equation, which is part of the quasiconformal theory. In this theory, Teichm\"uller space is realized as a quotient space of the unit ball in a Banach space of measurable Beltrami differentials over a fixed Riemann surface. The origin of the idea is in Teichm\"uller's 1939 paper \cite{T1939} although, as we said, Teichm\"uller was not aware of thie fact. For a concise survey on this point of view on the complex structure of Teichm\"uller space, we refer the reader to the paper \cite{E} by Earle.} The present paper is very different in spirit than the previous ones, and it did not attract much attention, compared to his previous papers on the subject, which were thoroughly analyzed and commented by Ahlfors, Bers and the schools they founded.\footnote{One reason why this paper did not attract much attention may be that Teichm\"uller's papers were examined by analysts and not by geometers, and the approach in this paper is algebro-geometric. Another reason is that the paper appeared in \emph{Deutsche Mathematik}, an ephemeral journal founded by Bieberbach to which very few libraries outside Germany subscribed. The journal contained, in the first two issues, articles presenting the Nazi viewpoint on the influence of race on mathematics. About the difficulty of reading Teichm\"uller's papers, we can quote Ahlfors from his 1954 paper \emph{On quasiconformal mappings} \cite{Ahlfors1954}, referring to the results of Teichm\"uller's 1939 paper \cite{T1939}: ``For the sake of completeness we have not hesitated to reproduce some of Teichm\"uller's reasonings almost without change. One good reason for this is that the Teichm\"uller papers are not easily available. Another reason is that it requires considerabe effort to extricate Teichm\"uller's complete and incontestable proofs from the maze of conjectures in which they are hidden."} The new approach to the moduli problem that Teichm\"uller presents in the paper that is the subject of this commentary is so different from the previous ones that is seems that Teichm\"uller himself was not sure that the space he obtains through the present techniques is the same as the space he introduced in his previous papers through the quasiconformal map approach.
That the paper was not read carefully by other mathematicians is testified by the fact that the paper is very rarely cited in the mathematical literature, and that there are results and methods in this paper that were rediscovered  later on without always referring to Teichm\"uller. Among these, we mention:
\begin{enumerate}
\item The existence and uniqueness of the universal Teichm\"uller curve, rediscovered later on by Ahlfors and by Bers. At the same time, this introduced the first fibre bundle over Teichm\"uller space.
\item The proof of the fact that the automorphisms group of the univeral Teichm\"uller curve is the extended mapping class group.
\item The idea of a fine moduli space.\footnote{A fine moduli space is a fiber space where the isomorphism type of the fibre determines the point below it.}
\item The idea of using the period map to define a complex structure on Teichm\"uller space.
\end{enumerate}
 
While commenting more on these facts, we shall give a quick review of the content of the paper.
 
Teichm\"uller starts by recalling that it was known before him, but only by heuristic arguments, that the number of complex parameters for the set of equivalence classes of Riemann surfaces of genus $g$ is 

$$\tau=\begin{cases} 0 & \hbox{ if } g=0 \\
              1 & \hbox{ if }  g=1 \\
              3(g-1) & \hbox{ if } g>1.
              \end{cases}
$$ 
He declares that several authors came up with these numbers using different methods, but that in reality these authors were not capable of saying precisely what they were counting. He considers that the fact that the various counts lead to the same value is a sort of a miracle, since the methods that were used were not rigorous. Teichm\"uller then says that these numbers, in order to be meaningful, should represent a dimension, and that in order to discuss the dimension of a set (in the present case, of moduli space), one has to turn this set into some ``space with a notion of neighborhood".\footnote{Teichm\"uller writes that ``as long as $\mathfrak{R}$ is not a space with a notion of neighborhood,
it does not have a dimension in the sense of analysis or set theory". In this context, the expression ``set theory"  means topology.} It appears here that Teichm\"uller was the first person who formulated in such precise terms the moduli problem for Riemann surfaces.

Teichm\"uller then emphasizes that one should not primarily ask for an explicit representation of points in the moduli space via numbers in a coordinate system (an approach which seems to have been suggested by Riemann's work), but that one should rather study the inner structure of that space. He then goes on saying that not only one would like to have on the set of moduli $\mathfrak{R}$ the structure of a topological space or of an algebraic variety, but one would also like to have the structure of an analytic manifold, that is to say, of a nonsingular complex space. 
Teichm\"uller says that this is not possible because $\mathfrak{R}$ ``contains certain singular manifolds". In modern language, this corresponds to the fact that moduli space is an orbifold and not a manifold. He therefore constructs a covering $\underline{\mathfrak{R}}$ of $\mathfrak{R}$ that has no singularities.  The space $\underline{\mathfrak{R}}$ is the space that was called later on \emph{Teichm\"uller space}.\footnote{It was probably Andr\'e Weil who first proposed the name ``Teichm\"uller space"; see Weil's comments in his  \emph{Collected Papers} (\cite{Weil-Collected} Vol. II, p. 546), where Weil writes: ``[...] this led me to the decision of writing up my observations on the moduli of curves and on what I called `Teichm\"uller space'  [...]"; see also the historical comments in \cite{JP}.}

Teichm\"uller then gives a short overview of his results and methods, and he notes that he will not be able to publish details ``in the near future".\footnote{Teichm\"uller died the same year, at the age of 30.} He announces that his solution to the problem of moduli is based on three newly introduced notions:
\begin{enumerate}
\item  \emph{The ``topological determination" of Riemann surfaces}:  this is the notion that we call today a ``marking" of a Riemann surface. Here, a Riemann surface is equipped with a fixed homotopy class of homeomorphisms from a fixed Riemann surface. We note that Teichm\"uller had already introduced markings in his 1939 paper \cite{T1939}. 
\item  \emph{The notion of an analytic family of Riemann surfaces}: this notion plays a central role in this paper as well as in  the later  developments of Teichm\"uller theory. It has been reintroduced later on by several authors working on moduli spaces, and below, we shall mention in particular Grothendieck.
\item \emph{The notion of ``turning piece coordinates"}: this is a operation of modifying the complex structure in the regular neighborhood of a simple closed curve in a Riemann surface. The complex automorphism group of the annulus is $S^1$ and the turning piece deformation can be considered in some sense as an ancestor of a complex analogue of the Fenchel-Nielsen deformation.
\end{enumerate}

After this introduction, Teichm\"uller defines the moduli space $\underline{\mathfrak{R}}$ and the action of the mapping class group on that space. He formulates the \emph{problem of moduli} as, a priori, the ``problem of asking for the properties of the space $\mathfrak{R}$". He says that however, it turns out that it is better to study, rather than the space $\mathfrak{R}$, its covering  $\underline{\mathfrak{R}}$. 
A formulation of the problem of moduli is again given in a more precise form.

Teichm\"uller then introduces the notion of an analytic $n$-dimensional manifold, defined by coordinate charts with holomorphic coordinate changes. It seems that this is one of the first appearances of such a definition in the mathematical  literature. We can quote here Remmert  (\cite{Remmert1998} p. 225):
\begin{quote} \small
It seems difficult to locate the first paper where complex manifolds explicitly occur. In 1944 they appear in Teichm\"uller's work on ``Ver\"anderliche Riemannsche Fl\"achen" (\emph{Collected Papers}, p. 714); here we find for the first time the German expression ``komplex analytische Mannigfaltigkeit". The English ``complex manifold" occurs in Chern's work (\cite{Chern1946} p. 103); he recalls the definition (by an atlas) just in passing. And in 1947 we find ``vari\'et\'e analytique complexe" in the title of Weil's paper \cite{Weil1947}. Overnight complex manifolds blossomed everywhere.
\end{quote}

In any case, it is an interesting fact that the first example of a complex manifold of higher dimension (other than the example of a domain of $\mathbb{C}^n$, $n\geq 2$) that appeared in the literature is precisely a space of (equivalence classes) of marked complex manifolds of dimension one.\footnote{Of course, Teichm\"uller space was not known yet to be a complex domain. An embedding of that space in a $\mathbb{C}^N$ was discovered later on.}

Teichm\"uller then introduces the notion of an \emph{analytic family of Riemann surfaces}. In modern language, this is a fiber bundle $\mathfrak{M}$ over an analytic base space $\mathfrak{B}$, the fibers being Riemann surfaces. The fiber bundle is locally trivial from the differentiable point of view (but not from the analytic point of view, since in a trivializing product neighborhood, two fibres are generally not isomorphic as Riemann surfaces).  

A particularly interesting analytic family of Riemann surfaces is the one where the base space $\mathfrak{B}$  is Teichm\"uller space and where the fibre above each point is a marked Riemann surface representing the point itself. In this case the fiber bundle is called, in modern language, the \emph{universal Teichm\"uller curve}, or the  \emph{Teichm\"uller curve}.

We note that Teichm\"uller curve has been re-introduced later on in the mathematical literature, in general with no reference to Teichm\"uller's paper.  

The fibre bundle approach to Teichm\"uller space was expanded and made precise by several authors, see e.g. Ahlfors \cite{Ahlfors-Some}, Bers \cite{Bers1958}, \cite{Grothendieck}  and Earle and Eells \cite{EE}. In Ahlfors' paper 1961 \cite{Ahlfors-Some} and in  Bers' 1961 paper \cite{Bers1958} this fiber bundle is used to define the complex structure on Teichm\"uller space.\footnote{Ahlfors, in his paper, says about his method (p. 171)  that ``this approach is essentially due to Bers", and Bers writes (p. 356)  that ``the existence of a `natural' complex structure in $\mathcal{T}_g$ has been asserted by Teichm\"uller; the first proof was given by Ahlfors after Rauch showed how to introduce complex-analytic co-ordinates in the neighborhood of any point which is not a hyperelliptic surface. Other proofs are due to Kodaira-Spencer and to Weil". Bers writes, at the beginning of the section on the analytic structure: ``The results of this and the following section confirm and extend some of Teichm\"uller's assertions in the paper \emph{Ver\"anderliche Riemannsche Fl\"achen}. They also show that the complex-analytic structure defined above is natural and coincides with that of Rauch-Ahlfors." We note by the way that the tangent bundle sequence of the bundle of Riemann surfaces is at the basis of the  Kodaira-Spencer  theory of infinitesimal deformations of Riemann surfaces  \cite{KS}, which also provides a description of the complex structure of Teichm\"uller space. Grothendieck's work on Teichm\"uller space is based on the consideration of analytic fiber bundles over surfaces (Grothendieck calls them ``algebraic curves of genus $g$), whose fibers are Riemann surfaces of genus $g$, cf. \cite{Grothendieck}, Expo\'e 7 and ff. Finally, let us note that the fiber bundle approach was also used to define analytic vector bundles, the most natural one being obtained by taking above each point of the base space the tangent bundle of the surface in the fiber. The study of characteristic classes of fiber bundles over moduli space, led to several important developments. We mention on this subject the Mumford conjecture (solved by Madsen) and the Witten conjecture (solved by Kontsevich).} 
The Teichm\"uller curve turned out to be an extremely important object in Hodge theory.

Teichm\"uller's aim was to show that the Teichm\"uller curve is a complex manifold of dimension $3g-2$.

In the general setting where the base  $\mathfrak{B}$ is an analytic manifold with fibers being Riemann surfaces, Teichm\"uller introduces a notion he calls \emph{permanent uniformizing local parameter},\footnote{This term, an English translation of the German ``permanente Ortsuniformisierende" is also used by Bers in \cite{Bers1958}.} which is a local analytic coordinate system that gives an analytic parameter on each of the fiber Riemann surfaces. Thus, Teichm\"uller gets a coordinate $t$ that works locally for a family of Riemann surfaces that are above points in $\mathfrak{B}$, which, together with the local $r$-dimensional parameters of $\mathfrak{B}$ produces a system of $(r+1)$-dimensional parameters of $\mathfrak{M}$ as an $(r+1)$-complex manifold. He states that such permanent uniformizing parameters exist.

The surfaces that are the fibers of the bundle over the space  $\mathfrak{B}$ are \emph{a priori} not marked. Teichm\"uller shows that from a marking on one fiber one can obtain a marking on nearby fibers. He then says that by well-known principles, from the space $\mathfrak{B}$ one can construct  a ``relatively unramified"  covering  $\underline{\mathfrak{B}}$ of  $\mathfrak{B}$  and where the surfaces above points of $\underline{\mathfrak{B}}$  are marked surfaces. 

Teichm\"uller then states an existence and uniqueness theorem for a globally analytic family of marked (Teichm\"uller says ``topologically determined") surfaces $\underline{\mathfrak{H}}[\frak{c}]$, where $\mathfrak{c}$ runs over a $\tau$-dimensional complex analytic manifold $\mathfrak{C}$ such that for any marked Riemann surface $\underline{\mathfrak{H}}$ of genus $g$ there is one and only one $\mathfrak{c}$ such that the Riemann surface  $\underline{\mathfrak{H}}$ is conformally equivalent to an $\underline{\mathfrak{H}}[\frak{c}]$ and such that the family $\underline{\mathfrak{H}}[\frak{c}]$ satisfies the following universal property: If $\underline{\mathfrak{H}}[\frak{p}]$ is any globally analytic family of Riemann surfaces with base $\mathfrak{B}$, there is a holomorphic map $f:\mathfrak{B}\to\mathfrak{C}$ such that the family $\underline{\mathfrak{H}}[\frak{p}]$ is the pull-back by $f$ of the family $\underline{\mathfrak{H}}[\frak{c}]$. Teichm\"uller states that such a family $\underline{\mathfrak{H}}[\frak{c}]$ exists and that it is \emph{essentially} unique. From the context, and stated in modern terms, essential uniqueness means that the family is unique up  to the action of the mapping class group.

The complex analytic manifold $\mathfrak{C}$, which is the base space of the family $\underline{\mathfrak{H}}[\frak{c}]$, is the object that we call today Teichm\"uller space. This existence and uniqueness result was rediscovered by Grothendieck, who gave a complete proof of it in an algebro-geometric language that  is different from Teichm\"uller's. Grothendieck gave a series of talks on this subject  at Cartan's seminar (1960-1961), and written texts of these talks were circulated and published. Grothendieck's statement is more general than that of Teichm\"uller; it is expressed in terms of a universal property, concerning (using Grothendieck's wording) a ``rigidifying functor" $\mathcal{P}$ relative to a discrete group $\gamma$, which can be taken in particular as the ``Teichm\"uller rigidifying functor", and where $\gamma$ is the mapping class group. Grothendieck's statement is the following (in this statement, $T$ is Teichm\"uller space.):
\begin{quote}\small
Theorem 3.1.---There exists an analytic space $T$ and a $\mathcal{P}$-algebraic curve $V$ above $T$ which are universal in the following sense: For every $\mathcal{P}$-algebraic curve $X$ above an analytic space $S$, there exists a unique morphism $g$ from $S$ to $T$ such that $X$ (together with its $\mathcal{P}$-structure) is isomorphic to the pull-back of $V/T$ by $g$.\footnote{(\cite{Grothendieck}  p. 7-08) [Th\'eor\`eme 3.1.--- Il existe un espace analytique $T$, et une $\mathcal{P}$-courbe alg\'ebrique $V$ au-dessus de $T$, qui soient universels au sens suivant : Pour toute $\mathcal{P}$-courbe alg\'ebrique $X$ au-dessus d'un espace analytique $S$, il existe un morphisme et un seul $g$ de $S$ dans $T$, tel que $X$ soit isomorphe (avec sa $\mathcal{P}$-structure) \`a l'image inverse par $g$ de $V/T$]. It seems that Grothendieck have heard of Teichm\"uller's papers, but like many others, he did not read them. Teichm\"uller's work on the problem of moduli had nevertheless an enormous influence on Grotendieck, who declares, in the introduction to this paper (Expos\'e 7, whose title is: ``An axiomatic description of Teichm\"uller space and its variants" [Description axiomatique de l'espace de Teichm\"uller et de ses variantes]: ``In doing this, the necessity of rewriting the foundations of analytic geometry will become manifest". [Chemin faisant,  la n\'ecessit\'e deviendra manifeste de revoir les fondements de la G\'eom\'etrie analytique]. It seems that the analysts working on Teichm\"uller theory had heard about Grothendieck's work, but did not understand it. We can quote here Abikoff, in a report he published in 1989 in the Bulletin of the AMS \cite{Abikoff}, on an book by Nag: ``First, algebraic geometers took us, the noble but isolated practitioners of this iconoclastic discipline, under their mighty wings. We learned the joys of providing lemmas solving partial differential and integral equations and various other nuts and bolts results. These served to render provable such theorems as: The ?\%$\sharp$\$! is representable." Let us also quote Ahlfors, from his 1964 survey on quasiconformal mappings \cite{A1964} (p. 152), talking about Teichm\"uller's 1944 paper: ``In a final effort Teichm\"uller produced a solution of the structure problem, by an entirely different method, but it was so cumbersome that it is doubtful whether anybody else has checked all the details. [...] It is only fair to mention, at this point, that the algebraists have also solved the problem of moduli, in some sense even more completely than the analysts. Because of the different language, it is at present difficult to compare the algebraic and analytic methods, but it would seem that both have their own advantages.}
\end{quote}

Grothendieck deduces the following corollary (in which $\gamma$ is, as before, the mapping class group):
\begin{quote}\small
Proposition 3.3.---  Let $X,X'$ be two $\mathcal{P}$-curves above $S$, defined respectively by morphisms $f,f'$ from $S$ in $T$. Assume that $S$ is connected and nonempty. Then the set of $S$-isomorphisms $X \stackrel{\sim}\to X'$ for the underlying curves (without the $\mathcal{P}$-structures), is in canonical one-to-one correspondence with the set of $u\in\gamma$ satisfying $f'=\overline{u}\circ f$.\footnote{(\cite{Grothendieck}  p. 7-10) [Proposition 3.3.--- Soient $X,X'$ deux $\mathcal{P}$-courbes au-dessus de $S$, d\'efinies respectivement par des morphismes $f,f'$ de $S$ dans $T$. Supposons $S$ connexe non vide. Alors l'ensemble des $S$-isomorphismes $X \stackrel{\sim}\to X'$ pour les courbes sous-jacentes (sans $\mathcal{P}$-structures), est en correspondance biunivoque canonique avec l'ensemble des $u\in\gamma$ tels que $f'=\overline{u}\circ f$].}
\end{quote}
This is a rigidity statement concerning the mapping class group. The mapping class group is canonically identified as the group of isomorphism classes of $\mathcal{P}$-curves.

This statement about the existence and uniqueness of the Teichm\"uller curve up to the mapping class group action can be considered as the first among a series of results that were obtained later on on the rigidity of mapping class group actions, the next (and probably the most famous) one being Royden's result stating that the automorphism group of the complex structure of Teichm\"uller space is the extended mapping class group, cf. \cite{Royden1971}. The existence and uniqueness of the universal family was later on constructed independently by Ahlfors and Bers, see the historical remarks in \cite{EM}.

Teichm\"uller considers that this theorem solves the \emph{moduli problem},\footnote{It is interesting to note that Grothendieck, after stating his main theorem (Theorem 3.1, stated above), makes a remark similar to that of Teichm\"uller (p. 7-10 of \cite{Grothendieck}): ``We shall see in the next sections that the analytic space $T$, equipped with the automorphism group $\gamma$, can be considered as a satisfying solution to the `moduli problem' for curves of genus $g$" [Nous verrons dans les paragraphes suivants que l'espace analytique $T$, muni du groupe d'automorphismes $\gamma$, peut \^etre consid\'er\'e comme une solution satisfaisant du  ``probl\`eme des modules" pour les courbes de genre $g$].} and on this occasion he formulates more precisely this problem: The space $\underline{\mathfrak{R}}$ of all classes of analytic marked surfaces of genus $g$ is made into a complex analytic manifold by identifying it to the base space $\mathfrak{C}$ of the  universal analytic family. This endows Teichm\"uller space at the same time with a complex analytic structure and with a topology.\footnote{Let us note that it was indeed considered, although this was not clearly stated, that the problem of moduli consists in defining a complex structure on moduli space. We can quote here Ahlfors, from his 1960 paper \cite{Ahlfors1960}, in which he writes the following: ``The classical problem [of moduli] calls for a complex analytic structure rather than a metric [...] The problem is not a clear cut one, and several formulations seem equally reasonable".} Teichm\"uller's proof of this theorem uses tools from algebraic geometry, as well as the algebro-geometric language of \emph{divisors}, \emph{principal parts} and \emph{places}. He states a second theorem which is needed in the proof of the first one,  result concerning  the ``determination of a function by generalized systems of principal parts".  In this contex, Riemann surfaces are studied through function fields. A Riemann surface is viewed as a field with degree of transcendence 1 over $\mathbb{C}$. The places in the field give back the surface.

 Teichm\"uller then introduces the notion of \emph{turning piece} and the one of \emph{turning piece coordinate} of a turning piece. This is a method of varying the complex structure of a surface which is given, in the tradition of Riemann, as a branched covering of the plane associated to an algebraic function. The analytic functions are determined by a ``generalized system of principal parts". The analytic structure of Teichm\"uller space is defined by varying the coefficients of these functions. 
 
 Teichm\"uller states that his space $\underline{\mathfrak{R}}$ ``consists of at most countably many connected parts", and he ``thinks that $\underline{\mathfrak{R}}$ in fact is simply connected". Thus, he was not sure that the space $\underline{\mathfrak{R}}$ he defines in this paper is the same as the Teichm\"uller space he defined in his 1939 paper \cite{T1939} using the quasiconformal theory.  It is interesting to note here that Grothendieck solved this issue. Indeed, after introducing the definition of Teichm\"uller space using the universal property, as we recalled above, Grothendieck writes, in (\cite{Grothendieck} 7-08):
 \begin{quote} \small
 It is also easy to check, using if needed a paper by Bers \cite{Bers-Tata}, that the space we introduce axiomatically here (if this space exists, and we shall prove this fact) is isomorphic to the Teichm\"uller space of the analysts. It follows that  Teichm\"uller space is homeomorphic to a ball, and therefore contractible, in particular connected and simply connected. A fortiori, the Jacobi spaces of all levels are connected, as is the moduli space $M$ introduced in Section 5 as a quotient space of Teichm\"uller space. It seems that at the time being there is no algebro-geometric proof even of the fact that moduli space is connected (which we can interpret in algebraic geometry by saying that two curves of the same genus $g$ are part of a family of algebraic curves parametrized by a connected algebraic variety).\footnote{[Il est d'ailleurs facile, de v\'erifier, utilisant au besoin un expos\'e de Bers \cite{Bers-Tata}, que l'espace introduit axiomatiquement ici (s'il existe, ce que nous prouverons) est isomorphe \`a l'espace de Teichm\"uller des analystes. Il en r\'esulte que l'espace de Teichm\"uller est hom\'eomorphe \`a une boule, et par suite contractile, en particulier connexe et simplement connexe. A fortiori, les espaces de Jacobi de tout \'echelon sont connexes, ainsi que l'espace des modules $M$ introduit au paragraphe 5 comme un espace quotient de l'espace de Teichm\"uller.  Il semble qu'il n'existe pas \`a l'heure actuelle de d\'emonstration, par voie alg\'ebrico-g\'eom\'etrique, m\^eme du fait que l'espace des modules est connexe, (ce qui s'interpr\`ete en g\'eom\'etrie alg\'ebrique en disant que deux courbes alg\'ebriques de m\^eme genre $g$ font partie d'une famille de courbes alg\'ebriques param\'etr\'ee par une vari\'et\'e alg\'ebrique connexe)].} 
 \end{quote}
  
 Then Teichm\"uller goes on studying the moduli space, that is, the space obtained by forgetting the marking. The description is very brief. (He says: ``At last I want to mention briefly what happens when releasing the topological determination".) He recalls the definition of  the mapping class group together with its action on the space of equivalence classes of marked surfaces, he shows that this action is properly discontinuous, and he defines the moduli space to be the quotient of Teichm\"uller space by this action. He considers what he calls the \emph{exceptional points} of $\underline{\mathfrak{R}}$; they lie on a certain analytic submanifold of $\underline{\mathfrak{R}}$. These points are stabilized by nontrivial finite groups of the mapping class groups, and they are obstructions for the moduli space to be a complex manifold. He shows that some of these points are \emph{substantial singularities}, that is, moduli space is not a manifold at these points. 
 
 In the last part of the paper, Teichm\"uller introduces the idea that the period matrices of differentials of the first kind can be used to study the complex structure of moduli space, and he states that hyperelliptic points may cause problems. This program on defining  the complex structure via the period map was carried on later by several authors including Rauch, who defined a complex structure away from hyperelliptic Riemann surfaces and showed that moduli space has non-uniformizable singularities at the surfaces which admit non-trivial conformal automorphisms, see \cite{Rauch} and \cite{Rauch-Weierstrass}. Ahlfors  \cite{Ahlfors-Some} also used the period map to define the complex structure of Teichm\"uller space, hence completing the program of Rauch.
 
 \bigskip

 {\bf Acknowledgements.} The authors are grateful to Bill Abikoff, Cliff Earle and Bill Harvey for valuable comments during the preparation of this commentary.

 \end{document}